\documentclass[12pt,reqno]{amsart}
\def\draft{n}
\def\published{n}
\usepackage{
  graphicx, amsmath, amssymb, multicol, mathtools, stmaryrd, xcolor,
  mathabx, needspace, lgreek, sidecap, mathbbol, blindtext, picins
}
\usepackage{enumitem}
\setlist[itemize]{left=0pt .. \parindent}
\setlist[enumerate]{left=0pt .. 1.5\parindent}
\usepackage{utfsym}     % for the likes of \car=\usym{1F697}
\usepackage[all]{xy}
\usepackage[textwidth=6.5in,textheight=9in,headsep=0.15in,centering]{geometry}

\usepackage[pagebackref,breaklinks=true]{hyperref}\hypersetup{colorlinks,
  linkcolor={green!50!black},
  citecolor={green!50!black},
  urlcolor=blue
}

\usepackage{chngcntr}
\counterwithin{figure}{section}

\usepackage{caption}
\if\draft y
  \usepackage[color]{showkeys}
\fi

\newcommand\ifpub[2]{\if\published y #1 \else #2 \fi}

\theoremstyle{plain}
\newtheorem{theorem}{Theorem}

\newtheorem{proposition}[theorem]{Proposition}

\newtheorem{lemma}[theorem]{Lemma}

\theoremstyle{definition}

\newtheorem{example}[theorem]{Example}

\newlength{\standardunitlength}
\def\car{\reflectbox{\usym{1F697}}}
\def\rac{\usym{1F697}}

\def\qed{{\linebreak[1]\null\hfill\text{$\Box$}}}

\newlength{\globalparindent}
\def\arXiv#1{{\href{https://arxiv.org/abs/#1}{arXiv:\linebreak[0]#1}}}

\def\weburl{http://drorbn.net/SL2PO}

\def\bbe{{\mathbb e}}

\def\bbH{{\mathbb H}}

\def\bbQ{{\mathbb Q}}

\def\calR{{\mathcal R}}
\def\calU{{\mathcal U}}
\def\fraka{{\mathfrak a}}

\def\frakg{{\mathfrak g}}

\def\sleps{{sl_{2+}^\eps}}

\def\draftcut{\if\draft y \cleardoublepage \fi}

\definecolor{lightred}{RGB}{255, 217, 217}

\def\imagetop#1{\vtop{\null\hbox{#1}}}

\def\eps{\epsilon}

\def\ip{{i^+}} \def\jp{{j^+}} \def\kp{{k^+}}
\def\ipp{{i^{+\!+}}} \def\jpp{{j^{+\!+}}} \def\kpp{{k^{+\!+}}}

\def\face{\begin{picture}(0,0)%
\includegraphics{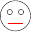}%
\end{picture}%
\setlength{\unitlength}{789sp}%
\begingroup\makeatletter\ifx\SetFigFont\undefined%
\gdef\SetFigFont#1#2#3#4#5{%
  \reset@font\fontsize{#1}{#2pt}%
  \fontfamily{#3}\fontseries{#4}\fontshape{#5}%
  \selectfont}%
\fi\endgroup%
\begin{picture}(1240,1240)(-19,-381)
\end{picture}%
}
\def\human{\begin{picture}(0,0)%
\includegraphics{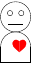}%
\end{picture}%
\setlength{\unitlength}{789sp}%
\begingroup\makeatletter\ifx\SetFigFont\undefined%
\gdef\SetFigFont#1#2#3#4#5{%
  \reset@font\fontsize{#1}{#2pt}%
  \fontfamily{#3}\fontseries{#4}\fontshape{#5}%
  \selectfont}%
\fi\endgroup%
\begin{picture}(1244,2442)(-21,-1583)
\end{picture}%
}
\def\machine{\begin{picture}(0,0)%
\includegraphics{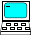}%
\end{picture}%
\setlength{\unitlength}{789sp}%
\begingroup\makeatletter\ifx\SetFigFont\undefined%
\gdef\SetFigFont#1#2#3#4#5{%
  \reset@font\fontsize{#1}{#2pt}%
  \fontfamily{#3}\fontseries{#4}\fontshape{#5}%
  \selectfont}%
\fi\endgroup%
\begin{picture}(1244,1394)(-21,-608)
\end{picture}%
}

\ifpub{\def\cellscale{0.76}}{\def\cellscale{1}}
\newsavebox\pdfbox
\newlength{\pdfheight}

\def\nbpdfInput#1{{%
  \savebox{\pdfbox}{\includegraphics[scale=\cellscale]{#1}}%
  \settoheight{\pdfheight}{\usebox{\pdfbox}}%
  \noindent\imagetop{\ifdim\pdfheight<10mm\begin{picture}(0,0)%
\includegraphics{figs/face.pdf}%
\end{picture}%
\setlength{\unitlength}{789sp}%
\begingroup\makeatletter\ifx\SetFigFont\undefined%
\gdef\SetFigFont#1#2#3#4#5{%
  \reset@font\fontsize{#1}{#2pt}%
  \fontfamily{#3}\fontseries{#4}\fontshape{#5}%
  \selectfont}%
\fi\endgroup%
\begin{picture}(1240,1240)(-19,-381)
\end{picture}%
\else\begin{picture}(0,0)%
\includegraphics{figs/human.pdf}%
\end{picture}%
\setlength{\unitlength}{789sp}%
\begingroup\makeatletter\ifx\SetFigFont\undefined%
\gdef\SetFigFont#1#2#3#4#5{%
  \reset@font\fontsize{#1}{#2pt}%
  \fontfamily{#3}\fontseries{#4}\fontshape{#5}%
  \selectfont}%
\fi\endgroup%
\begin{picture}(1244,2442)(-21,-1583)
\end{picture}%
\fi}\ %
  \imagetop{\usebox{\pdfbox}}%
  \vskip 2mm%
}}

\def\nbpdfEcho#1{{\noindent{\imagetop{\begin{picture}(0,0)%
\includegraphics{figs/machine.pdf}%
\end{picture}%
\setlength{\unitlength}{789sp}%
\begingroup\makeatletter\ifx\SetFigFont\undefined%
\gdef\SetFigFont#1#2#3#4#5{%
  \reset@font\fontsize{#1}{#2pt}%
  \fontfamily{#3}\fontseries{#4}\fontshape{#5}%
  \selectfont}%
\fi\endgroup%
\begin{picture}(1244,1394)(-21,-608)
\end{picture}%
}\ \imagetop{\includegraphics[scale=\cellscale]{#1}}\vskip 2mm}}}
\def\nbpdfMessage#1{{\noindent{\imagetop{\begin{picture}(0,0)%
\includegraphics{figs/machine.pdf}%
\end{picture}%
\setlength{\unitlength}{789sp}%
\begingroup\makeatletter\ifx\SetFigFont\undefined%
\gdef\SetFigFont#1#2#3#4#5{%
  \reset@font\fontsize{#1}{#2pt}%
  \fontfamily{#3}\fontseries{#4}\fontshape{#5}%
  \selectfont}%
\fi\endgroup%
\begin{picture}(1244,1394)(-21,-608)
\end{picture}%
}\ \imagetop{\includegraphics[scale=\cellscale]{#1}}\vskip 2mm}}}
\def\nbpdfOutput#1{{\noindent{\imagetop{\begin{picture}(0,0)%
\includegraphics{figs/machine.pdf}%
\end{picture}%
\setlength{\unitlength}{789sp}%
\begingroup\makeatletter\ifx\SetFigFont\undefined%
\gdef\SetFigFont#1#2#3#4#5{%
  \reset@font\fontsize{#1}{#2pt}%
  \fontfamily{#3}\fontseries{#4}\fontshape{#5}%
  \selectfont}%
\fi\endgroup%
\begin{picture}(1244,1394)(-21,-608)
\end{picture}%
}\ \imagetop{\includegraphics[scale=\cellscale]{#1}}\vskip 2mm}}}
\def\nbpdfPrint#1{{\noindent{\imagetop{\begin{picture}(0,0)%
\includegraphics{figs/machine.pdf}%
\end{picture}%
\setlength{\unitlength}{789sp}%
\begingroup\makeatletter\ifx\SetFigFont\undefined%
\gdef\SetFigFont#1#2#3#4#5{%
  \reset@font\fontsize{#1}{#2pt}%
  \fontfamily{#3}\fontseries{#4}\fontshape{#5}%
  \selectfont}%
\fi\endgroup%
\begin{picture}(1244,1394)(-21,-608)
\end{picture}%
}\ \imagetop{\includegraphics[scale=\cellscale]{#1}}\vskip 2mm}}}

\begin{document} %\latintext
\newdimen\captionwidth\captionwidth=\hsize
\setcounter{secnumdepth}{4}

\title{A Perturbed-Alexander Invariant}

\author{Dror~Bar-Natan}
\address{
  Department of Mathematics\\
  University of Toronto\\
  Toronto Ontario M5S 2E4\\
  Canada \ifpub{(corresponding author)}{}
}
\email{drorbn@math.toronto.edu}
\urladdr{http://www.math.toronto.edu/drorbn}

\author{Roland~van~der~Veen}
\address{
  University of Groningen, Bernoulli Institute\\
  P.O. Box 407\\
  9700 AK Groningen\\
  The Netherlands \ifpub{(corresponding author)}{}
}
\email{roland.mathematics@gmail.com}
\urladdr{http://www.rolandvdv.nl/}

\date{First edition June 24, 2022. This edition \today.}

\dedicatory{Dedicated to the memory of V.~F.~R.~Jones, 1952--2020, a friend and a mentor.}

\makeatletter
\@namedef{subjclassname@2020}{\textup{2020} Mathematics Subject Classification}
\makeatother

\thanks{This work was partially supported by NSERC grant RGPIN-2018-04350
and by the Chu Family Foundation (NYC). It is available in electronic
form, along with source files and a verification {\sl Mathematica}
notebook at \url{http://drorbn.net/APAI} and at \arXiv{2206.12298}.}

\begin{abstract}
In this note we give concise formulas, which lead to a simple and
fast computer program that computes a powerful knot invariant. This
invariant $\rho_1$ is not new, yet our formulas are by far the simplest
and fastest: given a knot we write one of the standard matrices $A$
whose determinant is its Alexander polynomial, yet instead of
computing the determinant we consider a certain quadratic expression in
the entries of $A^{-1}$. The proximity of our formulas to the Alexander
polynomial suggests that they should have a topological explanation. This
we do not have yet.

\end{abstract}

\maketitle

\setcounter{tocdepth}{3}
\tableofcontents
%\eject

\section{The Formulas} \label{sec:Formulas} One of the selling points for this
article is that the formulas in it are concise. Thus we start by running
through these formulas for a knot invariant $\rho_1$ as quickly as we
can. In Section~\ref{sec:Implementation} we turn the formulas into a
short yet very fast computer program, in Section~\ref{sec:Proofs} we
give a partial interpretation of the formulas in terms of car traffic
on a knot diagram and use it to prove the invariance of $\rho_1$, and
in Section~\ref{sec:Context} we quickly sketch the context: Alexander,
Burau, Jones, Melvin, Morton, Rozansky, Overbay, and our own prior
work. This article accompanies two talks,~\cite{Talk:Waco-2203} and~\cite{Talk:Geneva-2206}
(videos and handouts available).

\parpic[r]{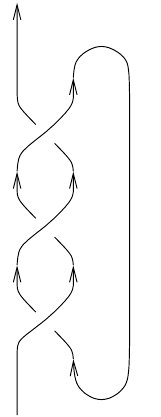}
Given an oriented $n$-crossing knot $K$, we draw it in the plane as a
long knot diagram $D$ in such a way that the two strands intersecting
at each crossing are pointed up (that's always possible because we
can always rotate crossings as needed), and so that at its beginning
and at its end the knot is oriented upward. We call such a diagram an
{\em upright knot diagram.} An example of an upright knot diagram $D_3$
is shown on the right.

We then label each edge of the diagram with two integer labels: a running index $k$ which runs from 1 to $2n+1$, and a ``rotation number'' $\varphi_k$, the geometric rotation number of that edge (the signed number of times the tangent to the edge is horizontal and heading right, with cups counted with $+1$ signs and caps with $-1$; this number is well defined because at their ends, all edges are headed up). On the right the running index runs from $1$ to $7$, and the rotation numbers for all edges are $0$ (and hence are omitted) except for $\varphi_4$, which is $-1$.

{\em A Technicality.} Some Reidemeister moves create or lose an edge and to avoid the need for renumbering it is beneficial to also allow labelling the edges with non-consecutive labels. Hence we allow that, and write $\ip$ for the successor of the label $i$ along the knot, and $\ipp$ for the successor of $\ip$ (these are $i+1$ and $i+2$ if the labelling is by consecutive integers). Also, ``$1$'' will always refer to the label of the first edge, and ``$2n+1$'' will always refer to the label of the last.

We let $A$ be the $(2n+1)\times(2n+1)$ matrix of Laurent polynomials in the formal variable $T$ defined by
\[ A \coloneqq I - \sum_c \left( T^sE_{i,\ip} + (1-T^s)E_{i,\jp} + E_{j,\jp} \right), \]
where $I$ is the identity matrix and $E_{\alpha\beta}$ denotes the elementary
matrix with $1$ in row $\alpha$ and column $\beta$ and zeros elsewhere.
The summation is over the crossings $c$ of the diagram $D$, and once $c$
is chosen, $s$ denotes its sign and $i$ and $j$ denote the labels below
the crossing where the label $i$ belongs to the over-strand and $j$
to the under-strand.

Alternatively, $A=I + \sum_c A_c$, where $A_c$ is a matrix of zeros except for the blocks as follows:
\begin{equation} \label{eq:A}
  \begin{array}{c}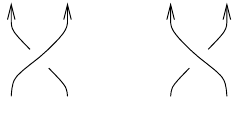\end{array}
  \qquad\longrightarrow\qquad
  \begin{array}{c|cccc}
    A_c &   \text{column }\ip  &  \text{column }\jp \\
    \hline
    \text{row }i & -T^s  & T^s-1 \\
    \text{row }j & 0  & -1
  \end{array}
\end{equation}

%{\color{blue} For example, if $D=D_1=\ \uparrow$ is the diagram with no crossings,
%the resulting matrix $A$ is the $1\times 1$ identity matrix
%$(1)$. \ldots}

\parpic[r]{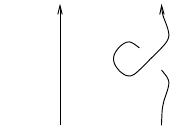}
For example, if $D=D_1$ is the diagram with no crossings (as shown on
the right), the resulting matrix $A$ is the $1\times 1$ identity matrix
$(1)$. If $D=D_2$ is the second diagram on the right (here $s=+1$,
$(i,j)=(2,1)$, and $(\ip,\jp)=(3,2)$), then
\[ A =
  \begin{pmatrix}1&0&0\\0&1&0\\0&0&1\end{pmatrix}
  + \begin{pmatrix}0&-1&0\\0&T-1&-T\\0&0&0\end{pmatrix}
  = \begin{pmatrix}1&-1&0\\0&T&-T\\0&0&1\end{pmatrix},
\]
and for $D_3$ as on the first page, we have
\[ A = \footnotesize \left(
  \begin{array}{ccccccc}
   1 & -T & 0 & 0 & T-1 & 0 & 0 \\
   0 & 1 & -1 & 0 & 0 & 0 & 0 \\
   0 & 0 & 1 & -T & 0 & 0 & T-1 \\
   0 & 0 & 0 & 1 & -1 & 0 & 0 \\
   0 & 0 & T-1 & 0 & 1 & -T & 0 \\
   0 & 0 & 0 & 0 & 0 & 1 & -1 \\
   0 & 0 & 0 & 0 & 0 & 0 & 1 \\
  \end{array}
\right). \]

\picskip{0}
We note without supplying details that the matrix $A$ comes in a
straightforward way from Fox calculus as it is applied to the Wirtinger
presentation of the fundamental group of the complement of $K$ (using
the diagram $D$). Hence the determinant of $A$ is equal up to a unit
to the normalized Alexander polynomial $\Delta$ of $K$ (which satisfies
$\Delta(T)=\Delta(T^{-1})$ and $\Delta(1)=1$). In fact, we have that
\begin{equation} \label{eq:Delta} \Delta = T^{(-\varphi(D)-w(D))/2}\det(A), \end{equation}
where $\varphi(D)\coloneqq\sum_k\varphi_k$ is the total rotation number of $D$ and where $w(D)=\sum_cs_c$ is the writhe of $D$, namely the sum of the signs $s_c$ of all the crossings $c$ in $D$.

For our example $D_2$, $\det(A)=T$, $\varphi(D)=1$, and $w(D)=1$, so $\Delta=T^{(-1-1)/2}\cdot T=1$, as expected for a diagram of the unknot.
For $D_3$, $\det(A)=1-T+T^2$, $\varphi(D)=-1$, and $w(D)=3$, so $\Delta=T^{(1-3)/2}(1-T+T^2)=T-1+T^{-1}$, as expected for the trefoil knot.

We set\footnote{At $T=1$ the matrix $A$ has $1$'s on the main diagonal, $(-1)$'s on the diagonal above it, and $0$'s everywhere else. Hence $A$ is invertible at $T=1$ and hence over the field of rational functions.} $G=(g_{\alpha\beta})=A^{-1}$, and taking our inspiration from physics, we name $g_{\alpha\beta}$ the {\em Green function} for the diagram $D$. For our three examples $D_1$, $D_2$, and $D_3$, the Green function $G$ is respectively
\begin{equation} \label{eq:GExamples}
  \begin{pmatrix} 1 \end{pmatrix},
  \quad \begin{pmatrix} 1&T^{-1}&1\\0&T^{-1}&1\\0&0&1 \end{pmatrix},
  \quad \footnotesize \left(
    \begin{array}{ccccccc}
     1 & \frac{T^3-T^2+T}{T^2-T+1} & 1 &
       \frac{T^3-T^2+T}{T^2-T+1} & 1 &
       \frac{T^3-T^2+T}{T^2-T+1} & 1 \\
     0 & 1 & \frac{1}{T^2-T+1} & \frac{T}{T^2-T+1} &
       \frac{T}{T^2-T+1} & \frac{T^2}{T^2-T+1} & 1 \\
     0 & 0 & \frac{1}{T^2-T+1} & \frac{T}{T^2-T+1} &
       \frac{T}{T^2-T+1} & \frac{T^2}{T^2-T+1} & 1 \\
     0 & 0 & \frac{1-T}{T^2-T+1} & \frac{1}{T^2-T+1} &
       \frac{1}{T^2-T+1} & \frac{T}{T^2-T+1} & 1 \\
     0 & 0 & \frac{1-T}{T^2-T+1} & \frac{T-T^2}{T^2-T+1}
       & \frac{1}{T^2-T+1} & \frac{T}{T^2-T+1} & 1 \\
     0 & 0 & 0 & 0 & 0 & 1 & 1 \\
     0 & 0 & 0 & 0 & 0 & 0 & 1 \\
    \end{array}
  \right).
\end{equation}

We can now define our invariant $\rho_1$. It is the sum of two sums. The first is a sum of a term $R_1(c)$ over all crossings $c$ in $D$, where for such a crossing we let $s$ denote its sign and we let $i$ and $j$ denote the edge labels of the incoming over- and under-strands, respectively and where
\begin{equation} \label{eq:R1}
  R_1(c) \coloneqq s \left(
    g_{ji} \left(g_{\jp,j}+g_{j,\jp}-g_{ij}\right)
    -g_{ii} \left(g_{j,\jp}-1\right)
    -1/2
  \right).
\end{equation}
The second sum is a sum over the edges $k$ of $D$ of a correction term dependent on the rotation number $\varphi_k$. We multiply the result by $\Delta^2$ to ``clear the denominators''\footnote{$R_1(s)$ is quadratic in the entries of $G$ and hence it has denominators proportional to $\Delta^2$.}:
\begin{equation} \label{eq:rho1}
  \rho_1 \coloneqq \Delta^2\left(\sum_c R_1(c) - \sum_k\varphi_k\left(g_{kk}-1/2\right)\right),
\end{equation}

Direct calculations show that $\rho_1(D_1)=0$ (as the sums are empty), $\rho_1(D_2)=0$, and $\rho_1(D_3) = -T^2 +2T - 2 + 2T^{-1} - T^{-2}$.

\ifpub{}{\needspace{12mm}}
\begin{theorem}[``Invariance'', proofs in Section~\ref{sec:Proofs}] \label{thm:Main} The quantity $\rho_1$ is a knot invariant.
\end{theorem}

As we shall see in the next section, $\rho_1$ has more separation power than the Jones polynomial, yet it is closer to the more topologically meaningful Alexander polynomial $\Delta$: it is cooked up from the same matrix $A$ and in terms of computational complexity, computing $\rho_1$ is not very different from computing $\Delta$. In order to compute $\Delta$ we need to compute the determinant of $A$, while to compute $\rho_1$ we need to invert $A$ and then compute a sum of $O(n)$ terms that are quadratic in the entries of $A^{-1}$.\footnote{We prefer not to be more specific about the complexity of computing $\rho_1$. It is the same as the complexity of inverting $A$, and matrix inversion is poly-time, with a rather small exponent, even for matrices with entries in a ring of polynomials (e.g.~\cite{Storjohann:ComplexityOfInversion}). We have not explored how much one can further gain by exploiting the fact that $A$ is very sparse.} We have computed $\rho_1$ for knots with over 200 crossings using the unsophisticated implementation presented in Section~\ref{sec:Implementation}.

Topologists should be intrigued! $\rho_1$ is derived from the same matrix as the Alexander polynomial $\Delta$, yet we have no topological interpretation for $\rho_1$.

\section{Implementation and Power} \label{sec:Implementation}

Two of the main reasons we like $\rho_1$ is that it is very easy to implement and even an unsophisticated implementation runs very fast. To highlight these points we include a full implementation here, a step-by-step run-through, and a demo run. We write in Mathematica~\cite{Wolfram:Mathematica}, and you can find the notebook displayed here at~\cite[APAI.nb]{Self}.

We start by loading the library \verb$KnotTheory`$~\cite{Bar-NatanMorrison:KnotTheory} (it is used here only for the list of knots that it contains, and to compute other invariants for comparisons). We also load a minor conversion routine~\cite[Rot.nb / Rot.m]{Self} whose internal workings are irrelevant here.

\noindent\nbpdfInput{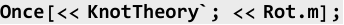}

\noindent\nbpdfPrint{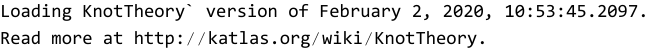}

\noindent\nbpdfPrint{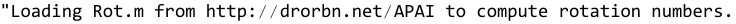}

\ifpub{}{\needspace{50mm}}
\subsection{The Program} This done, here is the full $\rho_1$ program:

\noindent\nbpdfInput{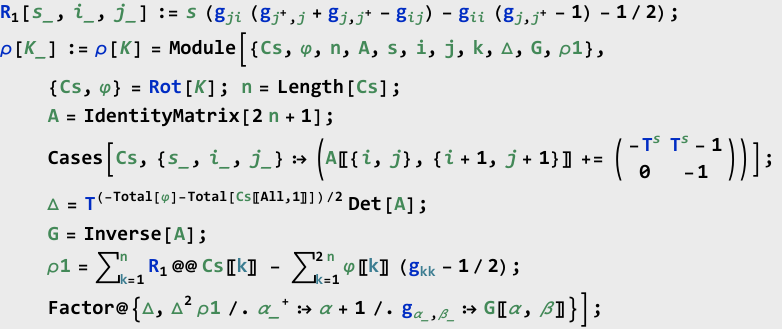}

The program uses mostly the same symbols as the text, so even without any knowledge of Mathematica, the reader should be able to recognize at least formulas~\eqref{eq:A}, \eqref{eq:Delta}, and~\eqref{eq:rho1} within it. As a further hint we add that the variable \verb$Cs$ ends up storing the list of crossings in a knot $K$, where each crossing is stored as a triple $(s,i,j)$, where $s$, $i$, and $j$ have the same meaning as in~\eqref{eq:A}. The conversion routine \verb$Rot$ automatically produces \verb$Cs$, as well as a list $\varphi$ of rotation numbers, given any other knot presentation known to the package \verb$KnotTheory`$.

Note that the program outputs the ordered pair $(\Delta,\rho_1)$. The Alexander polynomial $\Delta$ is anyway computed internally, and we consider the aggregate $(\Delta,\rho_1)$ as more interesting than any of its pieces by itself.

\subsection{A Step-by-Step Run-Through}  We start by setting $K$ to be the knot diagram on page~1 using the \verb$PD$ notation of \verb$KnotTheory`$~\cite{Bar-NatanMorrison:KnotTheory}. We then print \verb$Rot[K]$, which is a list of crossings followed by a list of rotation numbers:

\noindent\nbpdfInput{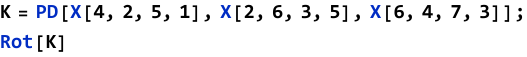}

\noindent\nbpdfOutput{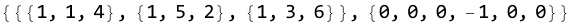}

Next we set \verb$Cs$ and $\varphi$ to be the list of crossings, and the list of rotation numbers, respectively.
\ifpub{}{\needspace{20mm}}

\noindent\nbpdfInput{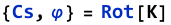}

\noindent\nbpdfOutput{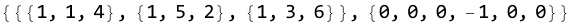}

We set \verb$n$ to be the number of crossings, \verb$A$ to be the $(2n+1)$-dimensional identity matrix, and then we iterate over \verb$c$ in \verb$Cs$,  adding a block as in~\eqref{eq:A} for each crossing.

\noindent\nbpdfInput{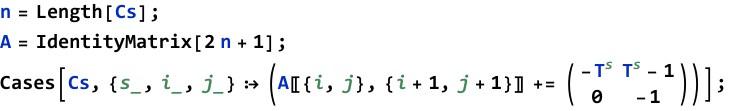}

\ifpub{}{\needspace{30mm}}
Here's what \verb$A$ comes out to be:

\noindent\nbpdfInput{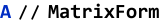}

\noindent\nbpdfOutput{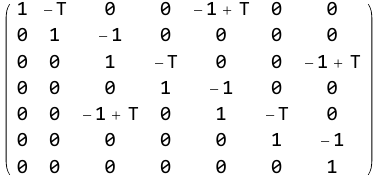}

We set $\Delta$ to be the determinant of \verb$A$, with a correction as in~\eqref{eq:Delta}. So $\Delta$ is the Alexander polynomial of $K$.

\noindent\nbpdfInput{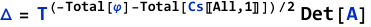}

\noindent\nbpdfOutput{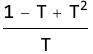}

\ifpub{}{\needspace{30mm}}
\verb$G$ is now the \verb$Inverse$ of \verb$A$:

\noindent\nbpdfInput{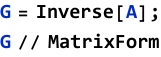}

\noindent\nbpdfOutput{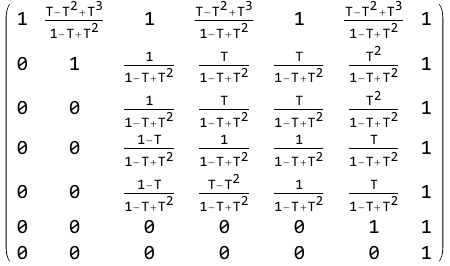}

\ifpub{}{\needspace{30mm}}
It remains to blindly follow the two parts of Equation~\eqref{eq:rho1}:

\noindent\nbpdfInput{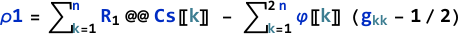}

\noindent\nbpdfOutput{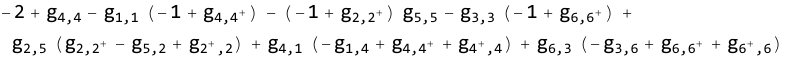}

We replace each ${\mathtt g}_{\alpha\beta}$ with the appropriate entry of \verb$G$:

\noindent\nbpdfInput{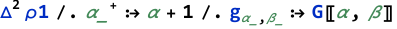}

\noindent\nbpdfOutput{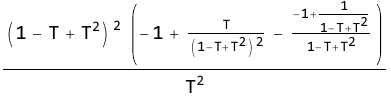}

Finally, we output both $\Delta$ and $\rho_1$. We factor them just to put them in a nicer form:

\noindent\nbpdfInput{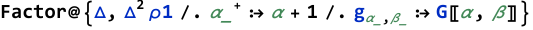}

\noindent\nbpdfOutput{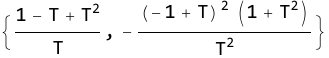}

\subsection{A Demo Run}  \label{ssec:Demo} Here are $\Delta$ and $\rho_1$ of all the knots with up to 6 crossings (a table up to 10 crossings is printed at~\cite{PG}:

\noindent\nbpdfInput{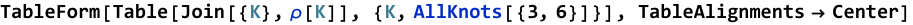}

\noindent\nbpdfMessage{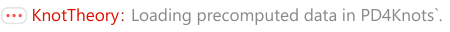}

\noindent\nbpdfOutput{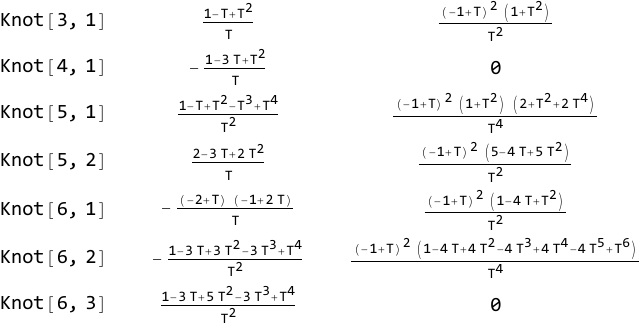}

Some comments are in order:
\begin{itemize}
\item If $\bar{K}$ is the mirror of a knot $K$, then $\rho_1(\bar{K})(T) = -\rho_1(K)(T^{-1})$. Indeed in~\eqref{eq:rho1} both $R_1(c)$ and $\varphi_k$ flip sign under reflection in a plane perpendicular to the plane of the knot diagram, and the matrix $A$ and hence also all the $g_{\alpha\beta}$'s are the same except for the substitution $T\to T^{-1}$.
\item $\rho_1$ seems to always be divisible by $(T-1)^2$ and seems to always be palindromic ($\rho_1(T)=\rho_1(T^{-1})$). We are not sure why this is so.
\item The last properties taken together would imply that $\rho_1$ vanishes on amphicheiral knots, such as $4_1$ and $6_3$ above.
\end{itemize}

\begin{figure}
\[ \resizebox{\ifpub{\linewidth}{6in}}{!}{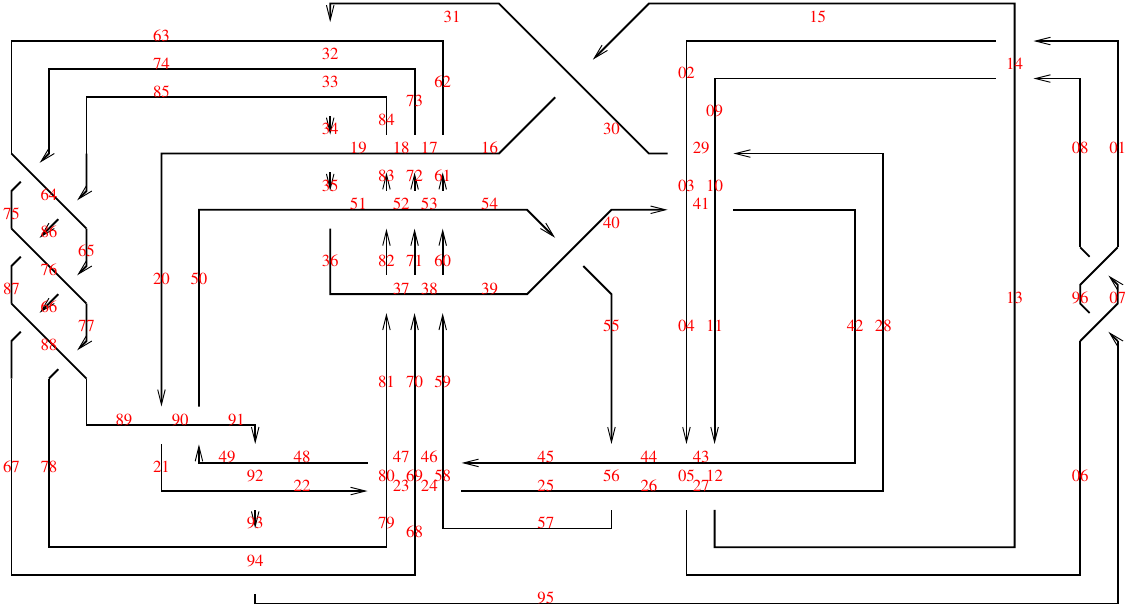} \]
\caption{A 48-crossing knot from~\cite{GompfScharlemannThompson:Counterexample}.} \label{fig:GST48}
\end{figure}
Next is one of our favourites, a knot from~\cite{GompfScharlemannThompson:Counterexample} (see Figure~\ref{fig:GST48}), which is a potential counterexample to the ribbon$=$slice conjecture~\cite{Fox:Problems}. It takes about two minutes to compute $\rho_1$ for this 48 crossing knot (note that Mathematica prints \verb$Timing$ information is seconds, and that this information is highly dependent on the CPU used, how loaded it is, and even on its temperature at the time of the computation):

\noindent\nbpdfInput{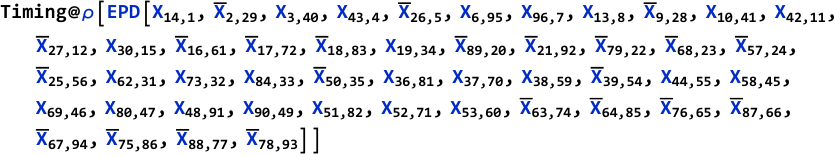}

\noindent\nbpdfOutput{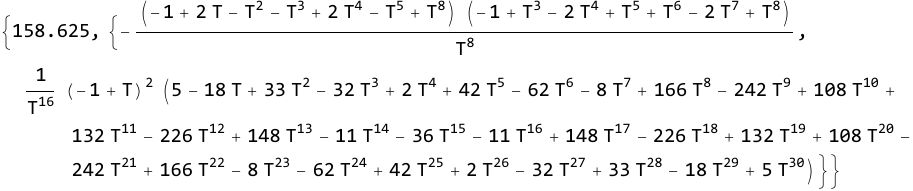}

\subsection{The Separation Power of $\rho_1$} Let us check how powerful $\rho_1$ is on knots with up to 12 crossings:

\noindent\nbpdfInput{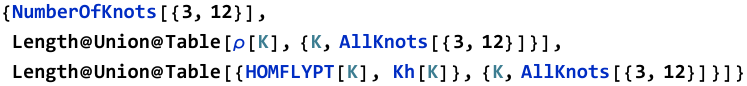}

\noindent\nbpdfOutput{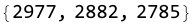}

So the pair $(\Delta,\rho_1)$ attains 2,882 distinct values on the 2,977 prime knots with up to 12 crossings (a deficit of 95), whereas the pair $(H,Kh)=$ (HOMFLYPT polynomial, Khovanov Homology) attains only 2,785 distinct values on the same knots (a deficit of 192).

In our spare time we computed all of these invariants on all the prime knots with up to 14 crossings. On these 59,937 knots the pair $(\Delta,\rho_1)$ attains 53,684 distinct values (a deficit of 6,253) whereas the pair $(H,Kh)$ attains only 49,149 distinct values on the same knots (a deficit of 10,788).

Hence the pair $(\Delta,\rho_1)$, computable in polynomial time by simple programs, seems stronger than the pair $(H,Kh)$, which is more difficult to program and (for all we know) cannot be computed in polynomial time. We are not aware of another poly-time invariant as strong as the pair $(\Delta,\rho_1)$.

\section{Proofs of Theorem~\ref{thm:Main}, the Invariance Theorem} \label{sec:Proofs}

We tell the proof of the Invariance Theorem (Theorem~\ref{thm:Main}) in two ways: an elegant and intuitive though slightly lacking telling in Section~\ref{ssec:Cars}, and a complete though slightly dull telling in Section~\ref{ssec:dull}. But first, a few common elements.

\subsection{Common Elements}

Two upright knot diagrams are considered the same (as diagrams) if their underlying knot diagrams are the same, and if the respective rotation numbers of their edges are all the same. It is clear that $\rho_1$ is well defined on upright knot diagrams. To prove Theorem~\ref{thm:Main} we need to know what to prove. Namely, when do two upright knot diagrams represent the same knot? This is answered in the spirit of the classical Reidemeister theorem by the following:

\begin{figure}
\ifpub
  {\[ \resizebox{\linewidth}{!}{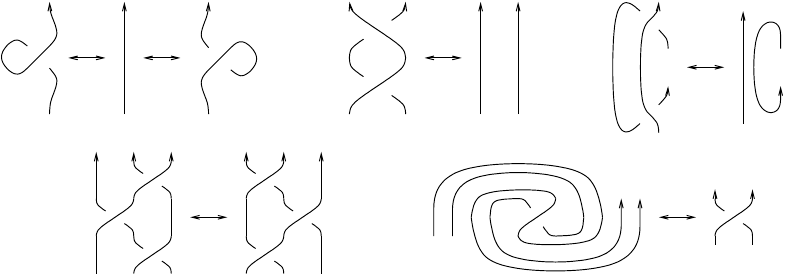} \]}
  {\[ \input{figs/UprightRMoves.pdf_t} \]}
\caption{
  The upright Reidemeister moves: Reidemeister 1 left and right,
  Reidemeister 2 braid-like and cyclic, Reidemeister 3, and (the $+$)
  Swirl.
} \label{fig:UprightRMoves}
\end{figure}

\begin{theorem}[``Upright Reidemeister''] \label{thm:UprightRMoves}
Two upright knot diagrams represent the same knot if and only if they
differ by a sequence of R1l, R2r, R2b, R2c, R3, and Sw$^+$ moves as in
Figure~\ref{fig:UprightRMoves}.
\end{theorem}

\noindent{\em Sketch of the proof.} In the case of round knots (i.e.,~
not ``long''), knot diagrams can be turned upright by rotating individual
crossings. The only ambiguity here is by powers of the full rotation,
the {\em swirls} Sw$^+$ and Sw$^-$ (where Sw$^-$ is the same as Sw$^+$
except with a negative crossing, and we don't need to impose it separately
as it follows from Sw$^+$ and R2). Hence we have a well-defined map
from $\{\text{knot diagrams}\}$ to $\{\text{upright knot diagrams}\}/$Sw$^\pm$.
It remains to write the usual Reidemeister moves between knot diagrams as
moves between upright knot diagrams. The result are the moves R1l, R2r,
R2b, R2c, and R3. Note that unoriented knot theory is presented with
just three Reidemeister moves, but these split into several versions
in the oriented case. The sufficiency of the versions we picked can be
found in~\cite{Polyak:RMoves}. In the case of long knots a minor further
complication arises, regarding the rotation numbers of the initial and
final edges. We leave the details of the problem and its resolution to
the reader. \qed

Our key formulas,~\eqref{eq:R1} and~\eqref{eq:rho1} involve the Green function $g_{\alpha\beta}$. We need to know that it is subject to some relations, the {\em $g$-rules} of Lemma~\ref{lem:gRules} below, whose proof is so easy that it comes first:

\noindent{\em Proof of Lemma~\ref{lem:gRules}.} The first set of $g$-rules
reads out column $\beta$ of the equality $AG=I$, and the second set of
$g$-rules reads out row $\alpha$ of the equality $GA=I$. \qed

\begin{lemma}[``$g$-rules''] \label{lem:gRules} Given a fixed upright knot diagram $D$, its corresponding matrix $A$, and its inverse $G=(g_{\alpha\beta})$, and given a crossing $c=(s,i,j)$ in $D$ (with $s$, $i$, and $j$ as before), the following two sets of relations (the $g$-rules) hold (with $\delta$ denoting the Kronecker delta):
\begin{equation} \label{eq:CarRules}
  g_{i\beta} = \delta_{i\beta}+T^sg_{\ip,\beta}+(1-T^s)g_{\jp,\beta},
  \qquad g_{j\beta} = \delta_{j\beta}+g_{\jp,\beta},
  \qquad g_{2n+1,\beta} = \delta_{2n+1,\beta}
\end{equation}
and
\begin{equation} \label{eq:CounterRules}
  g_{\alpha i} = T^{-s}(g_{\alpha,\ip}-\delta_{\alpha,\ip}),
  \qquad g_{\alpha j} = g_{\alpha,\jp} - (1-T^s)g_{\alpha i} - \delta_{\alpha,\jp},
  \qquad g_{\alpha,1} = \delta_{\alpha,1}.
\end{equation}
Furthermore, for each fixed $\beta$ there are $2n+1$ $g$-rules of type~\eqref{eq:CarRules} (the first two depend on a choice of one of $n$ crossings, and the third is fixed, to a total of $2n+1$ rules). These fully determine the $2n+1$ scalars $g_{\alpha\beta}$ corresponding to varying $\alpha$. Similarly, for each fixed $\alpha$ there are $2n+1$ $g$-rules of type~\eqref{eq:CounterRules}. These fully determine the $2n+1$ scalars $g_{\alpha\beta}$ corresponding to varying $\beta$.
\end{lemma}

For later use, we teach our computer about $g$-rules:

\noindent\nbpdfInput{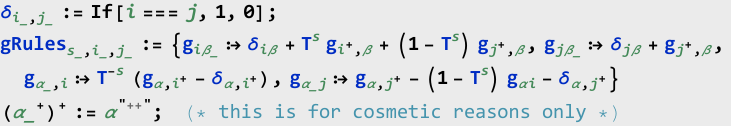}

\subsection{Cars, Traffic Counters, and Interchanges} \label{ssec:Cars}
Our first proof of Theorem~\ref{thm:Main} is slightly informal as it
uses the language and intuition of probability theory even though our
``probabilities'' are merely algebraic formulae and not numbers between
$0$ and $1$. Seasoned mathematicians should see that there is no real
problem here. Yet just to be safe, we also include a fully formal proof
in Section~\ref{ssec:dull}.

Cars ($\car$, $\rac$) travel on knot diagrams subject to the following three rules, inspired by Jones' ``bowling balls''~\cite{Jones:Hecke} and by Lin, Tian, and Wang's ``random walks''~\cite{LinTianWang:RandomWalk} (within the proof of Proposition~\ref{prop:CarsAreGreen} below we will see that these rules are equivalent to the $g$-rules of Equation~\eqref{eq:CarRules} above):
\begin{itemize}
\item On plain roads (edges) they travel following the orientation of the edge.
\item When reaching an underpass (the lower strand of a crossing), cars just pass through.
\item When reaching an overpass cars pass through with probability $T^s$
  (where $s=\pm 1$ is the sign of the crossing), yet drop over to the
  lower strand with the complementary probability of $1-T^s$.
\end{itemize}

These rules can be summarized by the following pictures:
\[ \resizebox{\linewidth}{!}{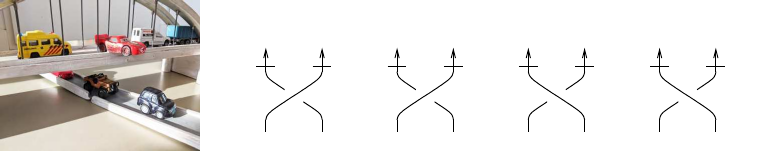} \]

In these pictures the horizontal struts represent ``traffic counters''
which measure the amount of traffic that passes through their respective
roads, and the output reading of these counters is printed above
them. Thus for example, the last interchange picture indicates that if
a unit stream of cars is injected into the diagram on the bottom right
and two traffic counters are placed at the top, then the first of these
will read a car intensity of $T^{-1}$ and the second $(1-T^{-1})$.

Note that our probabilities aren't really probabilities, if only because
$T$ and $T^{-1}$ cannot both be between $0$ and $1$ simultaneously. Yet
we will manipulate them algebraically as if they are probabilities,
restricting ourselves to equalities and avoiding inequalities. With
this restriction, we can use intuition from probability theory. We will
pretend that $T^s\sim 1$, or equivalently, that $1-T^s\sim 0$. This
has an algebraic meaning that does not refer to inequalities. Namely,
certain series can be deemed summable. For example,
\[ \sum_{r\geq 0}(1-T^s)^r = \frac{1}{1-(1-T^s)} = T^{-s}. \]

\parpic[r]{\input{figs/D2Traffic.pdf_T}}
\begin{example} Cars are injected on edge $\#1$ of the diagram $D_2$ of Section~\ref{sec:Formulas} as indicated on the right. What does the indicated traffic counter on edge $\#2$ measure?

\picskip{2}
\noindent{\em Solution.} Every car coming through the interchange from
$\#1$ passes through the underpass and comes to $\#2$, so the counter
reads ``$1$'' just for this traffic. But then these cars continue and
pass on the overpass, and $(1-T)$ of them fall down and continue through
edge $\#2$ and get counted again. But then these fallen cars continue and
pass on the overpass once again, and $(1-T)$ of them, meaning $(1-T)^2$
of the original traffic, fall once more and contribute a further reading
of $(1-T)^2$. This process continues and the overall counter reading is
\[ 1+(1-T)+(1-T)^2+(1-T)^3+\ldots = \frac{1}{1-(1-T)} = T^{-1}. \]
Note that this is exactly the row $1$ column $2$ entry of the matrix $G$ computed for this tangle in~\eqref{eq:GExamples}.
\end{example}

We claim that this is general:

\begin{proposition} \label{prop:CarsAreGreen} For a general knot diagram $D$, the entry
$g_{\alpha\beta}$ of its Green function is equal to the reading of a
traffic counter placed at $\beta$ given that traffic is injected into $D$
at $\alpha$. (In the case where $\alpha=\beta$, the counter is placed
after where the traffic is injected, not before).
\end{proposition}

\noindent{\em Proof.} Consider the $g$-rules of type~\eqref{eq:CarRules}. The third, $g_{2n+1,\beta} = \delta_{2n+1,\beta}$ is the statement that if traffic is injected on the outgoing edge of $D$, it can only be measured on the outgoing edge of $D$ (so traffic never flows backwards). The second, $g_{j\beta} = \delta_{j\beta}+g_{\jp,\beta}$, is the statement that traffic goes through underpasses undisturbed, so $g_{j\beta} = g_{\jp,\beta}$ unless the traffic counter $\beta$ is placed between $j$ and $\jp$, in which case it measures one unit more if the cars are injected before it, at $j$, rather than after it, at $\jp$. Similarly the first of these $g$-rules, $g_{i\beta} = \delta_{i\beta}+T^sg_{\ip,\beta}+(1-T^s)g_{\jp,\beta}$, is the statement of the behaviour of traffic at overpasses. Thus the rules in~\eqref{eq:CarRules} are obeyed by cars and traffic counters, and as the rules in~\eqref{eq:CarRules} determine $g_{\alpha\beta}$, the proposition follows. \qed

\begin{proposition} The quantity $\rho_1$ is invariant under R3.
\end{proposition}

\noindent{\em Proof.} We first show that cars entering a multiple
interchange styled as the left hand side of the R3 move, exit it with
the exact same distribution as cars entering the multiple interchange
styled as the right hand side. The hardest part of that computation is
when cars enter at the bottom left (at $i$) and it boils down to the
equality $1\!-\!T = (1\!-\!T)^2\!+\!T(1\!-\!T)$:
\[
  \def\messA{{$(1\!-\!T)^2\!+\!T(1\!-\!T)$}}
  \def\messB{{$(1\!-\!T)T$}}
  \def\messC{{$T(1\!-\!T)$}}
  \def\messD{{$1\!-\!T$}}
  \ifpub
    {\resizebox{0.95\linewidth}{!}{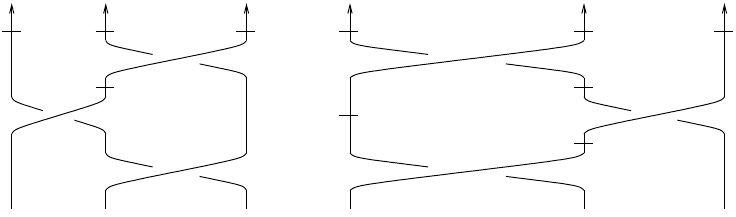}}
    {\input{figs/R3Traffic.pdf_t}}
\]
If cars enter in the middle or at the bottom right (at $j$ or at $k$), the computation is even easier.\footnote{Note that this computation is exactly the one that proves that the Burau representation~\cite{Burau:Zopfgruppen} respects R3.}

The conclusion is that performing the R3 move does not affect traffic patterns outside the area of the move itself; namely, the Green function $g_{\alpha\beta}$ is unchanged if both $\alpha$ and $\beta$ are outside the area of the move.

Thus the only contribution to $\rho_1$ that may change
(see~\eqref{eq:rho1}) is the contribution coming from the three $R_1$
terms corresponding to the crossings that move, and we need to know if
the following equality holds:
\ifpub
{\begin{multline*}
  R_1(+1,j,k) + R_1(+1,i,\kp) + R_1(+1,\ip,\jp) \\
  \overset{?}{=} R_1(+1,i,j) + R_1(+1,\ip,k) + R_1(+1,\jp,\kp)
\end{multline*}}
{\[
  R_1(+1,j,k) + R_1(+1,i,\kp) + R_1(+1,\ip,\jp)
  \overset{?}{=} R_1(+1,i,j) + R_1(+1,\ip,k) + R_1(+1,\jp,\kp)
\]}

Both sides here are messy quadratics involving the $g_{\alpha\beta}$'s of both sides, evaluated at $\alpha,\beta\in\{i,j,k,\ip,\jp,\kp,\ipp,\jpp,\kpp\}$. But we can use the traffic rules (aka the $g$-rules) to rewrite these quadratics in terms of the $g_{\alpha\beta}$'s with $\alpha,\beta\in\{\ipp,\jpp,\kpp\}$, and these are unchanged between the sides. So we simply need to know whether the above equality holds {\em after} the relevant $g$-rules have been applied to both sides. We could do that by hand, but it's simpler to appeal to a higher wisdom:

\noindent\nbpdfInput{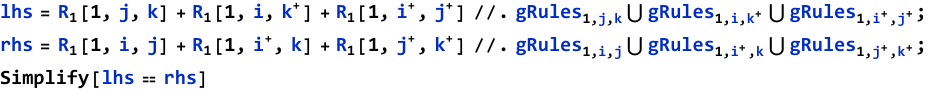}

\noindent\nbpdfOutput{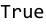}

 \ \qed

\noindent{\em First Proof of Theorem~\ref{thm:Main}, ``Invariance''.} We've shown invariance under R3. Invariance under the other moves is shown in a similar way: first one shows that overall traffic patterns are unchanged by each of the moves, and then one verifies that the local contributions to $\rho_1$ coming from the area changed by each move are equal once the $g$-rules are used to rewrite them in terms of $g_{\alpha\beta}$'s that are unaffected by the moves. This is shown in greater detail in the following section. \qed

\subsection{A More Formal Version of the Proof} \label{ssec:dull}
Again we start with the hardest, R3:

\begin{proposition} \label{prop:R3} The quantity $\rho_1$ is invariant under R3.
\end{proposition}

\noindent{\em Proof.} We need to know how the Green function
$g_{\alpha\beta}$ changes under R3. Here are the two sides of the move,
along with the $g$-rules of type~\eqref{eq:CarRules} corresponding to
the crossings within, written with the assumption that $\beta$ isn't in
$\{\ip,\jp,\kp\}$, so several of the Kronecker deltas can be ignored. We
use $g$ for the Green function at the left-hand side of R3, and $g'$ for the
right-hand side:
\[
  \def\grulesA{{$\begin{array}{l}
    g_{j,\beta} = \delta_{j\beta} \!+\! Tg_{\jp,\beta} \!+\! (1\!-\!T)g_{\kp,\beta} \\
    g_{k,\beta} = \delta_{k\beta} \!+\! g_{\kp,\beta}
  \end{array}$}}
  \def\grulesB{{$\begin{array}{l}
    g_{i,\beta} = \delta_{i\beta} \!+\! Tg_{\ip,\beta} \!+\! (1\!-\!T)g_{\kpp,\beta} \\
    g_{\kp,\beta} = g_{\kpp,\beta}
  \end{array}$}}
  \def\grulesC{{$\begin{array}{l}
    g_{\ip,\beta} = Tg_{\ipp,\beta} \!+\! (1\!-\!T)g_{\jpp,\beta} \\
    g_{\jp,\beta} = g_{\jpp,\beta}
  \end{array}$}}
  \def\gprulesA{{$\begin{array}{l}
    g'_{i,\beta} = \delta_{i\beta} \!+\! Tg'_{\ip,\beta} \!+\! (1\!-\!T)g'_{\jp,\beta} \\
    g'_{j,\beta} = \delta_{j\beta} \!+\! g'_{\jp,\beta}
  \end{array}$}}
  \def\gprulesB{{$\begin{array}{l}
    g'_{\ip,\beta} = Tg'_{\ipp,\beta} \!+\! (1\!-\!T)g'_{\kp,\beta} \\
    g'_{k,\beta} = \delta_{k\beta} \!+\! g'_{\kp,\beta}
  \end{array}$}}
  \def\gprulesC{{$\begin{array}{l}
    g'_{\jp,\beta} = Tg'_{\jpp,\beta} \!+\! (1\!-\!T)g'_{\kpp,\beta} \\
    g'_{\kp,\beta} = g'_{\kpp,\beta}
  \end{array}$}}
  \ifpub
    {\resizebox{\linewidth}{!}{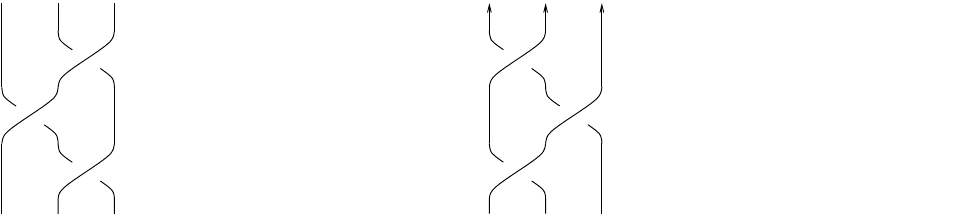}}
    {\input{figs/R3.pdf_t}}
\]

A routine computation (eliminating $g_{\ip,\beta}$, $g_{\jp,\beta}$, and $g_{\kp,\beta}$) shows that the first system of 6 equations is equivalent to the following 3 equations:
\[ g_{i,\beta} = \delta_{i\beta} + T^2g_{\ipp,\beta} + T(1-T)g_{\jpp,\beta} + (1-T)g_{\kpp,\beta}, \]
\[
  g_{j,\beta} = \delta_{j\beta} + Tg_{\jpp,\beta} + (1-T)g_{\kpp,\beta},
  \qquad\text{and}\qquad
  g_{k,\beta} = \delta_{k\beta} + g_{\kpp,\beta}.
\]
Similarly eliminating $g'_{\ip,\beta}$, $g'_{\jp,\beta}$, and $g'_{\kp,\beta}$ from the second set of equations, we find that it is equivalent to
\[ g'_{i,\beta} = \delta_{i\beta} + T^2g'_{\ipp,\beta} + T(1-T)g'_{\jpp,\beta} + (1-T)g'_{\kpp,\beta}, \]
\[
  g'_{j,\beta} = \delta_{j\beta} + Tg'_{\jpp,\beta} + (1-T)g'_{\kpp,\beta},
  \qquad\text{and}\qquad
  g'_{k,\beta} = \delta_{k\beta} + g'_{\kpp,\beta}.
\]
But these two sets of equations are the same, and as stated in the
$g$-rules lemma (Lemma~\ref{lem:gRules}), along with the $g$-rules
corresponding to the other crossings in $D$ (which are also the same
between $g$ and $g'$), these equations determine $g_{\alpha\beta}$
and $g'_{\alpha\beta}$, for $\alpha,\beta\notin\{\ip,\jp,\kp\}$. So
with this exclusion on $\alpha$ and $\beta$, we have that
$g_{\alpha\beta}=g'_{\alpha\beta}$. But this means that the
summations~\eqref{eq:rho1} in the definitions of $\rho_1$ are equal for
the two sides of R3, except perhaps for the three summands on each side
that come from the crossings that touch $\{\ip,\jp,\kp\}$.

What remains is completely mechanical. We just need to compute the sum of those three summands for both sides of R3, and apply to it the $g$-rules of types~\eqref{eq:CarRules} and~\eqref{eq:CounterRules} that eliminate the indices $\{\ip,\jp,\kp\}$. The computation is easy enough to be done by hand, yet why bother? Here's the machine version (it takes less typing to apply all relevant $g$-rules and also eliminate the indices $\{i,j,k\}$):

\ifpub{}{\needspace{30mm}}

\noindent\nbpdfInput{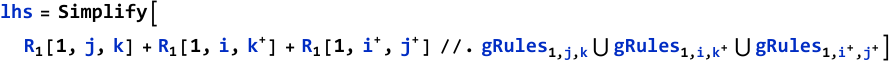}

\noindent\nbpdfOutput{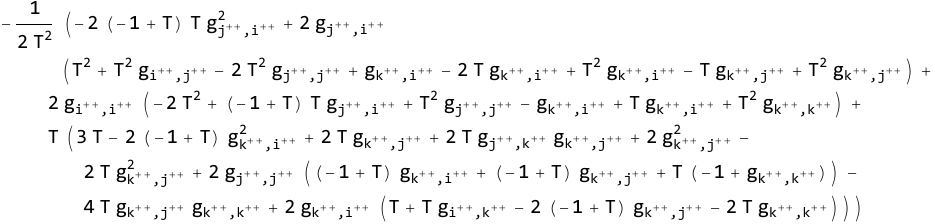}

\noindent\nbpdfInput{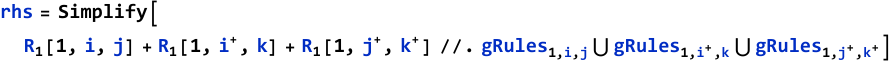}

\noindent\nbpdfOutput{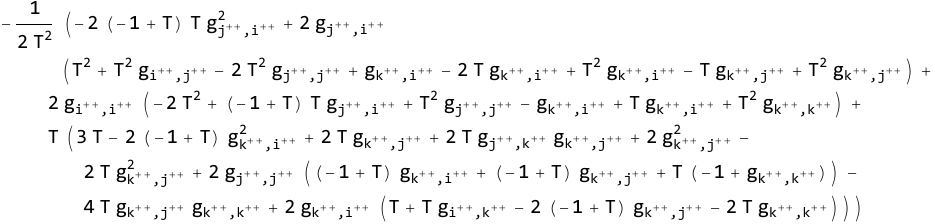}

\noindent\nbpdfInput{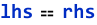}

\noindent\nbpdfOutput{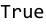}

 \ \qed

\begin{proposition} \label{prop:R2c} The quantity $\rho_1$ is invariant under R2c.
\end{proposition}

\noindent{\em Proof.} We follow the exact same steps as in the case
of $R3$. First, we write the $g$-rules, assuming that $\beta$ is not in
$\{i,j,\ip,\jp\}$:
\[
  \def\grulesA{{$\begin{array}{l}
    g_{i,\beta} = T^{-1}g_{\ip,\beta}+(1-T^{-1})g_{\jpp,\beta} \\
    g_{\jp,\beta} = g_{\jpp,\beta}
  \end{array}$}}
  \def\grulesB{{$\begin{array}{l}
    g_{\ip,\beta} = Tg_{\ipp,\beta}+(1-T)g_{\jp,\beta} \\
    g_{j,\beta} = g_{\jp,\beta}
  \end{array}$}}
  \def\rotA{{$\varphi_{j^+}\!=\!1$}}
  \def\rotB{{$\varphi_{j^{+\!+}}\!=\!1$}}
  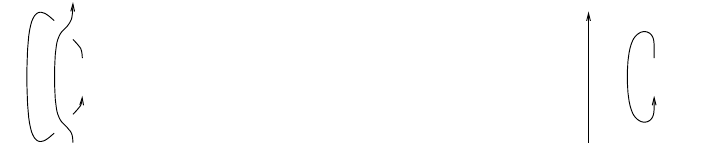
\]

Note that for the right hand side we allowed ourselves to label the edges $\ipp$ and $\jpp$ as the computation is independent of the labelling and the labelling need not be by contiguous integers (outside of the move area, we assume that the left hand side and the right hand side are labelled in the same way). Note also that for the right hand side, there are no relevant $g$-rules. Now as in the case of R3, for the left hand side we eliminate $g_{\ip,\beta}$ and $g_{\jp,\beta}$ and we are left with the relations
\[ g_{i,\beta}=g_{\ipp,\beta} \qquad\text{and}\qquad g_{j,\beta}=g_{\jpp,\beta}. \]
Otherwise the $g$-rules for the left and for the right are the same, and so their Green functions are the same except if the indices are in $\{i,j,\ip,\jp\}$ (these indices do not even appear in the right hand side). Thus the contribution to $\rho_1$ from outside the area of the move is the same for both sides.

Next we write the contribution to $\rho_1$ coming from the two crossings
and one rotation that appear on the left, and use the $g$-rules to push
all the indices in $\{i,j,\ip,\jp\}$ up to $\ipp$ and $\jpp$. This can
be done by hand, but seeing that we have tools, we use them as follows:

\noindent\nbpdfInput{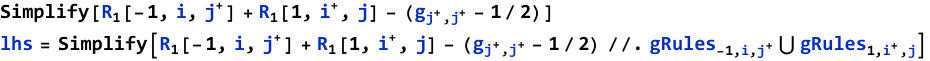}

\noindent\nbpdfOutput{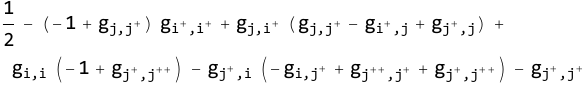}

\noindent\nbpdfOutput{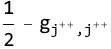}

This result is clearly equal to the single rotation contribution to $\rho_1$ that comes from the right hand side. \qed

\begin{proposition} \label{prop:R1l} The quantity $\rho_1$ is invariant under R1l.
\end{proposition}

\noindent{\em Proof.} We start with the relevant $g$-rules:
\[
  \def\grules{{$\begin{array}{l}
    g_{\ip,\beta} = \delta_{\ip,\beta} + Tg_{\ipp,\beta} + (1-T)g_{\ip,\beta} \\
    g_{i,\beta}=\delta_{i,\beta}+g_{\ip,\beta}
  \end{array}$}}
  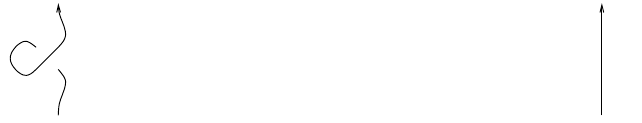
\]
The first of these rules is equivalent to $g_{\ip,\beta}=T^{-1}\delta_{\ip,\beta} + g_{\ipp,\beta}$. For $\beta\neq i,\ip$ we find as before that $g_{i,\beta}=g_{\ipp,\beta}$ and we can ignore the contributions to $\rho_1$ coming from outside the area of the move. The contribution to $\rho_1$ coming from the single crossing and single rotation on the left hand side is computed below, and is equal to the empty contribution coming from the right hand side:

\noindent\nbpdfInput{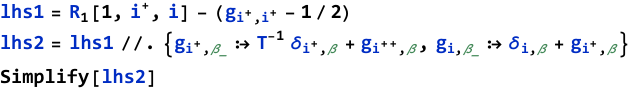}

\noindent\nbpdfOutput{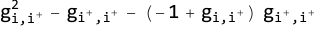}

\noindent\nbpdfOutput{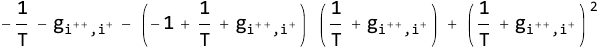}

\noindent\nbpdfOutput{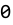}

 \ \qed

\vskip 1mm
\noindent{\em Second Proof of Theorem~\ref{thm:Main}, ``Invariance''.} After the Upright Reidemeister Theorem (Theorem~\ref{thm:UprightRMoves}) which sets out what we need to do, and Propositions~\ref{prop:R3}, \ref{prop:R2c}, and~\ref{prop:R1l} which prove invariance under R3, R2c, and R1l, it remains to show the invariance of $\rho_1$ under R1r, R2b, and Sw$^+$. This is done exactly as in the examples already shown, so in each case we show only the punch line:

\noindent\nbpdfInput{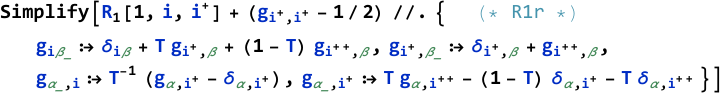}

\noindent\nbpdfOutput{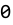}

\noindent(Note that the version of the $g$-rules we used above easily follows from~\eqref{eq:CounterRules}).

\noindent\nbpdfInput{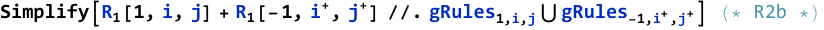}

\noindent\nbpdfOutput{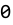}

\noindent\nbpdfInput{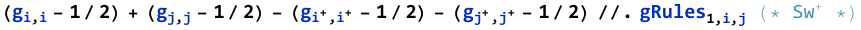}

\noindent\nbpdfOutput{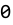}

 \ \qed

\section{Some Context and Some Morals} \label{sec:Context} We would
like to emphasize again that $\rho_1$ seems very close to the Alexander
polynomial, yet we have no topological interpretation for it. Until that
changes, where is $\rho_1$ coming from?

It comes via a lengthy path, which we will only sketch here. For
a while now \cite{PP1, PG, Talk:NCSU-1604, Talk:Indiana-1611,
Talk:LesDiablerets-1708, Talk:Matemale-1804, Talk:Toronto-1811,
Talk:DaNang-1905, Talk:UCLA-191101} we've been studying quantum invariants
related to the Lie algebra $\sleps$, the 4-dimensional Lie algebra with generators $y,b,a,x$ and
brackets
\[ [b,x]=\eps x, \quad [b,y]=-\eps y, \quad [b,a]=0, \quad [a,x]=x, \quad [a,y]=-y, \quad [x,y]=b+\eps a, \]
where $\eps$ is a scalar. The beauty of this algebra stems from the following:
\begin{itemize}
\item It is a ``classical double'' of a two-dimensional the Lie bialgebra $\langle a,x\rangle$, with
  \[ [a,x]=x,\qquad \delta(a)=0, \qquad \delta(x)=\eps x\wedge a, \]
  and hence quantization tools are available and are used below (e.g.~\cite{EtingofSchiffman:QuantumGroups}).
\item At invertible $\eps$ it is isomorphic to $sl_2\oplus\langle
  t\rangle$, where $t$ is a central element\footnotemark. Quantum
  topology tells us that the algebra $sl_2$ is related to the
  Jones polynomial. In fact, the universal quantum invariant
  (see \cite{Lawrence:UniversalUsingQG, Lawrence:Universal,
  Ohtsuki:QuantumInvariants}) for the Lie algebra $sl_2$ is equivalent
  to the coloured Jones polynomial of~\cite{Jones:Hecke}.
\item At $\eps=0$ it becomes the diamond Lie algebra $\diamond$, a solvable
  algebra in which computations are easier. The algebra $\diamond$
  is the semi-direct product of the unique non-commutative 2D Lie
  algebra $\fraka$ with its dual, and quantum topology tells us that it
  is related to the Alexander polynomial~\cite{Talk:Chicago-1009, WKO1}.
\end{itemize}
The last two facts taken together tell us that the Alexander
polynomial is some limit of the coloured Jones polynomial
(originally conjectured~\cite{MelvinMorton:Coloured, Rozansky:TrivialFlatI} and proven
by other means~\cite{Bar-NatanGaroufalidis:MMR}).

\footnotetext{{\arraycolsep=1pt\def\arraystretch{1}
Via the isomorphism
$\left(\begin{array}{cc}1&0\\0&-1\end{array}\right) \leftrightarrow \eps^{-1}b+a$,
$\left(\begin{array}{cc}0&1\\0&0\end{array}\right) \leftrightarrow x$,
$\left(\begin{array}{cc}0&0\\1&0\end{array}\right) \leftrightarrow \eps^{-1}y$,
and $t \leftrightarrow b-\eps a$.
}}

We can make this a bit more explicit.  By using the Drinfel'd quantum
double construction~\cite{Drinfeld:QuantumGroups} we find that the
universal enveloping algebra $\calU(\sleps)$ has a quantization $QU$,
which has an $R$-matrix solving the Yang-Baxter equation (meaning,
satisfying the R3 move, in the appropriate sense). These are given by:
\[ QU =
  A\langle y,b,a,x\rangle\left/\begin{pmatrix}
    [b,a]=0, \quad [b,x]=\eps x,\quad [b,y]=-\eps y, \\
    [a,x]=x, \quad [a,y]=-y, \quad xy-qyx=\frac{1-\bbe^{-\hbar(b+\eps a)}}{\hbar}
  \end{pmatrix}\right.,
\]
where $A\langle\text{\it gens}\rangle$ is the free associative algebra
with generators {\it gens}, and $q=\bbe^{\hbar\eps}$, and
\begin{multline*}
   R=\sum_{m,n\geq 0}\frac{y^nb^m\otimes (\hbar a)^m(\hbar x)^n}{m![n]_q!}
  \ifpub{\\}{\qquad}
  \left(\text{where $[n]_q!=\prod_{k=1}^n\frac{1-q^k}{1-q}$ is a ``quantum factorial''}\right).
\end{multline*}

Thus there is an associated universal quantum invariant of knots $Z_\eps(K)\in QU$
(which, as stated, is equivalent to the coloured Jones polynomial). In
our talks and papers we show that $Z_\eps$ can be expanded as a power
series in $\eps$, that at $\eps=0$ it is equivalent to the Alexander
polynomial, and that in general, the coefficient $Z_{(k)}$ of $\eps^k$
in $Z_\eps$ can be computed in polynomial time and is homomorphic,
meaning that it leads to an ``algebraic knot theory'' in the sense of
(say)~\cite{Talk:Sydney-190916}. We also know that the excess information in $Z_{(k)}$ (beyond the
information in $\{Z_{(0)},\ldots,Z_{(k-1)}\}$) is contained in a single polynomial, $\rho_k$. The first of
these polynomials is $\rho_1$ of this paper.

But how did we arrive at the specific formulas of this paper? As
often seen with quantizations, $QU$ is isomorphic (though only as
an algebra, not as a Hopf algebra) with $\calU(\sleps)$, and the
latter can be represented into the Heisenberg algebra
\ifpub
  {\[ \bbH=A\langle p,x\rangle/([p,x]=1) \]}
  {$\bbH=A\langle p,x\rangle/([p,x]=1)$}
via
\[ y\to -tp-\eps\cdot xp^2,\qquad b\to t+\eps\cdot xp,\qquad a\to xp,\qquad x\to x, \]
(abstractly, $\sleps$ acts on its Verma module
$\calU(\sleps)/(\calU(\sleps)\langle y,a,b-\eps a-t\rangle)\cong\bbQ[x]$
by differential operators, namely via $\bbH$). So $QU$'s $R$-matrix
can be expanded in powers of $\eps$ and pushed to $\calU(\sleps)$ and on
to $\bbH$, resulting in $\calR=\calR_0(1+\eps \calR_1+\cdots)$, with
$\calR_0=\bbe^{t(xp\otimes 1-x\otimes p)}$ and $\calR_1$ a quartic
polynomial in $p$ and $x$. Now all the computations for $\rho_1$ can be carried out by pushing around a
rather small number of $p$'s and $x$'s (at most 4), and this can be done using the rules
\begin{align*}
  (p\otimes 1)\calR_0 &= \calR_0(\bbe^t(p\otimes 1)+(1-\bbe^t)(1\otimes p)), \\
  (1\otimes p)\calR_0 &= \calR_0(1\otimes p),
\end{align*}
which, after setting $T=\bbe^t$, must remind the reader of
Equation~\eqref{eq:A}. When all the dust settles, the resulting formulas
are similar to the ones in Equations~\eqref{eq:R1} and~\eqref{eq:rho1}
(but only similar, because we applied some ad hoc cosmetics to make the
formulas appear nicer).

There are some morals to this story:

\begin{enumerate}

\item The definition of $\rho_1$ in Section~\ref{sec:Formulas} and the
proofs of its invariance in Section~\ref{sec:Proofs} are clearly much
simpler than the origin story, as outlined above. So quite clearly, we
still don't understand $\rho_1$. There ought to be a room for it directly
within topology, which does not require that one would know anything about
quantum algebra. (And better if that room is large enough to accommodate
morals~(\ref{mor:rhok}) and~(\ref{mor:frakg}) below).

\item \label{mor:rhok}
Like there is $\rho_1$, there are $\rho_k$. The origin story tells us that
$\rho_k$ should have a formula as a summation over choices of $k$-tuples
of features of the knot (crossings and rotations), just as the formula
for $\rho_1$ is a single summation over these features. The summand
for $\rho_k$ will be a degree $2k$ polynomial in the Green function
$g_{\alpha\beta}$ (compare with~\eqref{eq:R1}, which is quadratic). As a
$k$-fold summation, after inverting $A$, $\rho_k$ should be computable in
$O(n^k)$ additions and multiplications of polynomials in $T$, where $n$
is the crossing number.

\item These $\rho_k$ should be equivalent to the invariants in our earlier
works \cite{PP1, PG, Talk:NCSU-1604, Talk:Indiana-1611,
Talk:LesDiablerets-1708, Talk:Matemale-1804, Talk:Toronto-1811,
Talk:DaNang-1905, Talk:UCLA-191101}.

\item These $\rho_k$ should be equivalent to the invariants studied
earlier by Rozansky and Overbay \cite{Rozansky:TrivialFlatI,
Rozansky:Burau, Rozansky:U1RCC, Overbay:Thesis}, as their quantum
origin is essentially the same (though strictly speaking, we have not
written proofs of that, and normalizations may differ). Our formulas are
significantly simpler and faster to compute than the Rozansky-Overbay
formulas, and in our language it is easier to see the behaviour of
$\rho_1$ under mirror reflection (see Section~\ref{ssec:Demo}).

\item Like the Rozansky-Overbay invariants, $\rho_k$ should be equivalent to the ``higher diagonals'' for the
Melvin-Morton expansion (e.g.~\cite{Rozansky:U1RCC}) and should be dominated by the ``loop expansion'' of
the Kontsevich integral~\cite{Kricker:Lines, GaroufalidisRozansky:LoopExpansion}.

\item \label{mor:frakg} The quantum algebra story extends to other Lie
algebras, beyond $sl_2$. So there should be variants $\rho_k^\frakg$
of $\rho_k$ at least for every semisimple Lie algebra $\frakg$, given by
more or less similar formulas. Quantum algebra suggests that
$\rho_k^\frakg$ should be a polynomial in as many variables as the
rank of $\frakg$, and should in general be stronger than the ``base''
$\rho_k$. We have not seriously explored $\rho_k^\frakg$ yet, though
some preliminary work was done by Schaveling~\cite{Schaveling:Thesis}.

\item It appears that $QU$ has interesting traces and therefore that there should be a link version of
$\rho_1$. We have not pursued this formally.

\item $QU$ has a co-product and an antipode, and so the universal
tangle invariant associated with $QU$ has formulas for strand reversal
and strand doubling (e.g.~\cite{PG, Talk:DaNang-1905}). This implies
(e.g., by following the ideas of~\cite{Talk:Sydney-190916}) that there
should be formulas for $\rho_1$ that start with a Seifert surface for
the knot. We are pursuing such formulas now; we already know that the
degree of $\rho_1$ is bounded by $2g$, where $g$ is the genus of a knot~\cite{PG}.

\item For the same reasons, for ribbon knots $\rho_1$ should
have a formula computable from a ribbon presentation, and its
values might be restricted in a manner similar to the Fox-Milnor
condition~\cite{FoxMilnor:CobordismOfKnots}. We are pursuing this now.

\item The coloured Jones polynomial is invariant under mutation so we
expect $\rho_1$ to likewise be invariant under mutation (and indeed, also
$\rho_k$), yet we do not have a direct proof of that yet. Note that we can expect $\rho_k^\frakg$ for
higher-rank $\frakg$ to no longer be invariant under mutation.

\end{enumerate}


\begin{thebibliography}{LTW}
\backrefparscanfalse\renewcommand*\backref[1]{See pp.~#1.}%
\renewcommand*\backrefalt[4]{\ifcase #1\or In page #2.\else In pages #2.\fi}
\bibitem[BN1]{Talk:Chicago-1009} D.~Bar-Natan,
  {\em From the $ax+b$ Lie Algebra to the Alexander Polynomial and Beyond,}
  talk given at Knots in Chicago, September 2010. Video and handout at \url{http://drorbn.net/Chi10}.
  \backrefprint

\bibitem[BN2]{Talk:NCSU-1604} D.~Bar-Natan,
  {\em Gauss-Gassner Invariants,}
  talk given at Knots in the Triangle (Knots in Washington XLII), North Carolina State University, April 29
  -- May 1, 2016. Video and handout at \url{http://drorbn.net/NCSU-1604}. \backrefprint

\bibitem[BN3]{Talk:Indiana-1611} D.~Bar-Natan,
  {\em A Poly-Time Knot Polynomial Via Solvable Approximation,}
  talk given at Indiana University, November 7, 2016. Video and handout at
  \url{http://drorbn.net/Indiana-1611}. \backrefprint

\bibitem[BN4]{Talk:LesDiablerets-1708} D.~Bar-Natan,
  {\em The Dogma is Wrong,}
  talk given at Lie Groups in Mathematics and Physics, Les Diablerets, August 2017. Video and handout at
  \url{http://drorbn.net/ld17}. \backrefprint

\bibitem[BN5]{Talk:Toronto-1811} D.~Bar-Natan,
  {\em Computation without Representation,}
  talk given in Toronto, November 21, 2018. Video and handout at
  \url{http://drorbn.net/to18}. \backrefprint

\bibitem[BN6]{Talk:DaNang-1905} D.~Bar-Natan,
  {\em Everything around $\sleps$ is DoPeGDO. So what?,}
  talk given at Quantum Topology and Hyperbolic Geometry Conference, Da Nang, Vietnam, May 27--31 2019.
  Video and handout at \url{http://drorbn.net/v19}. \backrefprint

\bibitem[BN7]{Talk:Sydney-190916} D.~Bar-Natan,
  {\em Algebraic Knot Theory,}
  talk given in Sydney, September 16, 2019. Video and handout at \url{http://drorbn.net/syd2}.
  \backrefprint

\bibitem[BN8]{Talk:UCLA-191101} D.~Bar-Natan,
  {\em Some Feynman Diagrams in Algebra,}
  talk given in UCLA, November 1, 2019. Video and handout at \url{http://drorbn.net/la19}. \backrefprint

\bibitem[BN9]{Talk:Waco-2203} D.~Bar-Natan,
  {\em Cars, Interchanges, Traffic Counters, and a Pretty Darned Good Knot Invariant,}
  talk given at the 55th Spring Topology and Dynamical
  Systems Conference, Waco, March 2022. Video and handout at
  \url{http://drorbn.net/waco22}. \backrefprint

\bibitem[BN10]{Talk:Geneva-2206} D.~Bar-Natan,
  {\em Cars, Interchanges, Traffic Counters, and a Pretty Darned Good Knot Invariant,}
  talk given at ``From Subfactors to Quantum Topology - In memory of Vaughan Jones'',
  Geneva, June 30 2022. Video and handout at
  \url{http://drorbn.net/j22}. \backrefprint

\bibitem[BD]{WKO1} D.~Bar-Natan and Z.~Dancso,
  \href
    {http://drorbn.net/AcademicPensieve/Projects/WKO1}
    {\em Finite Type Invariants of W-Knotted Objects I: W-Knots and the Alexander Polynomial,}
  Alg.\ and Geom.\ Top.\ {\bf 16-2} (2016) 1063--1133,
  \arXiv{1405.1956}. \backrefprint

\bibitem[BNG]{Bar-NatanGaroufalidis:MMR} D.~Bar-Natan and S.~Garoufalidis,
  {\em On the Melvin-Morton-Rozansky conjecture,}
  Invent.\ Math.\ {\bf 125} (1996) 103--133. \backrefprint

\bibitem[BV1]{PP1} D.~Bar-Natan and R.~van~der~Veen,
  {\em A Polynomial Time Knot Polynomial,}
  Proc.\ Amer.\ Math.\ Soc.\ {\bf 147} (2019) 377--397, \arXiv{1708.04853}. \backrefprint

\bibitem[BV2]{Talk:Matemale-1804} D.~Bar-Natan and R.~van~der~Veen,
  {\em Talks in Matemale,}
  April 2018. Videos and handout at \url{http://drorbn.net/mm18}. \backrefprint

\bibitem[BV3]{PG} D.~Bar-Natan and R.~van~der~Veen,
  {\em Perturbed Gaussian Generating Functions for Universal Knot Invariants,}
  \arXiv{2109.02057}. \backrefprint

\bibitem[BV4]{Self} D.~Bar-Natan and R.~van der Veen,
  {\em A Perturbed-Alexander Invariant,}
  (self-reference), paper and related files at
  \url{http://drorbn.net/APAI}. The published and the \arXiv{2206.12298} editions may be older.
  \backrefprint

\bibitem[BM]{Bar-NatanMorrison:KnotTheory} D.~Bar-Natan, S.~Morrison, and others,
  {\tt KnotTheory`}, a knot theory mathematica package,
  \url{http://katlas.org/wiki/The_Mathematica_Package_KnotTheory}. \backrefprint

\bibitem[Bu]{Burau:Zopfgruppen} W.~Burau,
  {\em \"Uber Zopfgruppen und gleichsinnig verdrillte Verkettungen,}
  Abh.\ Math.\ Semin.\ Univ.\ Hamburg {\bf 11} (1936) 179–-186. \backrefprint

\bibitem[Dr]{Drinfeld:QuantumGroups} V.~G.~Drinfel'd,
  {\em Quantum Groups,}
  Proc.\ Int.\ Cong.\ Math., 798--820, Berkeley, 1986. \backrefprint

\bibitem[ES]{EtingofSchiffman:QuantumGroups} P.~Etingof and O.~Schiffman,
  {\em Lectures on Quantum Groups,}
  International Press, Boston, 1998. \backrefprint

\bibitem[Fo]{Fox:Problems} R.~H.~Fox,
  {\em Some Problems in Knot Theory,}
  Topology of 3-manifolds and related topics (Proc.\ The Univ. of Georgia Institute, 1961), Prentice-Hall, Englewood Cliffs, N.J., 168--176. \backrefprint

\bibitem[FM]{FoxMilnor:CobordismOfKnots} R.~H.~Fox and J.~W.~Milnor,
  {\em Singularities of 2-Spheres in 4-Space and Cobordism of Knots,}
  Osaka J.\ Math.\ {\bf 3-2} (1966) 257--267. \backrefprint

\bibitem[GR]{GaroufalidisRozansky:LoopExpansion} S.~Garoufalidis and L.~Rozansky,
  {\em The Loop Expansion of the Kontsevich Integral, the Null Move and $S$-Equivalence,}
  Topology {\bf 43-5} (2004) 1183--1210, \arXiv{math/0003187}. \backrefprint

\bibitem[GST]{GompfScharlemannThompson:Counterexample} R.~E.~Gompf, M.~Scharlemann, and A.~Thompson,
  {\em Fibered Knots and Potential Counterexamples to the Property 2R and Slice-Ribbon Conjectures,}
  Geom.\ and Top.\ {\bf 14} (2010) 2305--2347, \arXiv{1103.1601}. \backrefprint

\bibitem[Jo]{Jones:Hecke} V.~F.~R.~Jones,
  {\em Hecke Algebra Representations of Braid Groups and Link Polynomials,}
  Annals Math., {\bf 126} (1987) 335-388. \backrefprint

\bibitem[Kr]{Kricker:Lines} A.~Kricker,
  {\em The Lines of the Kontsevich Integral and Rozansky's Rationality Conjecture,}
  \arXiv{math/0005284}. \backrefprint

\bibitem[La1]{Lawrence:UniversalUsingQG} R.~J.~Lawrence,
  {\em Universal Link Invariants using Quantum Groups,}
  Proc\. XVII Int.\ Conf.\ on Diff.\ Geom.\ Methods in Theor.\ Phys., Chester,
  England, August 1988. World Scientific (1989) 55--63. \backrefprint

\bibitem[La2]{Lawrence:Universal} R.~J.~Lawrence,
  {\em A Universal Link Invariant,}
  Proc.\ IMA Conf.\ Math.\ --- Particle Phys.\ Interface, Oxford, England,
  September 1988. Inst.\ Math.\ Appl.\ Conf.\ Ser.\ New Ser.\ {\bf 24},
  Oxford Univ.\ Press (1990) 151--156. \backrefprint

\bibitem[LTW]{LinTianWang:RandomWalk} X-S.~Lin, F.~Tian, and Z.~Wang,
  {\em Burau Representation and Random Walk on String Links,}
  Pac.\ J.\ Math., {\bf 182-2} (1998) 289--302, \arXiv{q-alg/9605023}. \backrefprint

\bibitem[MM]{MelvinMorton:Coloured} P.~M.~Melvin and H.~R.~Morton,
  {\em The coloured Jones function,}
  Comm.\ Math.\ Phys.\ {\bf 169} (1995) 501--520. \backrefprint

\bibitem[Oh]{Ohtsuki:QuantumInvariants} T.~Ohtsuki,
  {\em Quantum Invariants,}
  Series on Knots and Everything {\bf 29}, World Scientific 2002. \backrefprint

\bibitem[Ov]{Overbay:Thesis} A.~Overbay,
  {\em Perturbative Expansion of the Colored Jones Polynomial,}
  Ph.D.\ thesis, University of North Carolina, August 2013,
  \url{https://cdr.lib.unc.edu/concern/dissertations/hm50ts889}. \backrefprint

\bibitem[Po]{Polyak:RMoves} M.~Polyak,
  {\em Minimal Generating Sets of Reidemeister Moves,}
  Quantum Topol.\ {\bf 1} (2010) 399--411, \arXiv{0908.3127}. \backrefprint

\bibitem[Ro1]{Rozansky:TrivialFlatI} L.~Rozansky,
  {\em A Contribution of the Trivial Flat Connection to the Jones
  Polynomial and Witten's Invariant of 3D Manifolds, I,}
  Comm.\ Math.\ Phys.\ {\bf 175-2} (1996) 275--296, \arXiv{hep-th/9401061}. \backrefprint

\bibitem[Ro2]{Rozansky:Burau} L.~Rozansky,
  {\em The Universal $R$-Matrix, Burau Representation and the Melvin-Morton
    Expansion of the Colored Jones Polynomial,}
  Adv.\ Math.\ {\bf 134-1} (1998) 1--31, \arXiv{q-alg/9604005}. \backrefprint

\bibitem[Ro3]{Rozansky:U1RCC} L.~Rozansky,
  {\em A Universal $U(1)$-RCC Invariant of Links and Rationality Conjecture,}
  \arXiv{math/0201139}. \backrefprint

\bibitem[Sch]{Schaveling:Thesis} S.~Schaveling,
  {\em Expansions of Quantum Group Invariants,}
  Ph.D.\ thesis, Universiteit Leiden, September 2020,
  \url{https://scholarlypublications.universiteitleiden.nl/handle/1887/136272}. \backrefprint

\bibitem[St]{Storjohann:ComplexityOfInversion} A.~Storjohann,
  {\em On the Complexity of Inverting Integer and Polynomial Matrices,}
  Comput.\ Complex.\ {\bf 24} (2015) 777-–821. \backrefprint

\bibitem[Wo]{Wolfram:Mathematica}
  {\em Wolfram Language \& System Documentation Center,}
  \url{https://reference.wolfram.com/language/}. \backrefprint

\end{thebibliography}
\end{document}
\endinput